\documentclass[reqno,tbtags,a4paper,12pt]{amsart}

\usepackage[in]{fullpage}

\usepackage{amsmath,amssymb,amsthm,amsfonts,amscd}

\newtheorem{theorem}{Theorem}
\newtheorem{lemma}{Lemma}

\renewcommand{\Bbb}[1]{\mathbb{#1}}

\author{A.\,A.\,Gaifullin}

\thanks{The work was partially supported by the Russian Foundation for Basic Research, grant
no.~05-01-01032, and the Russian Leading Scientific School
Support, grant no.~4182.2006.1.}

\title{Explicit construction of manifolds realizing the prescribed
homology classes}

\email{gaifull@mccme.ru}

\date{}

\address{Moscow State University}

\begin{document}
\maketitle

We consider a well-known N.\,Steenrod's problem on realization of
homology classes by images of the fundamental classes of
manifolds. In this paper we give an explicit combinatorial
construction that for a given cycle produces a piecewise linear
manifold realizing a multiple of the homology class of this cycle.
The construction is based on a local procedure for resolving
singularities of the cycle. An approach to N.\,Steenrod's problem
based on resolving singularities of pseudo-manifolds was
initially proposed by D.\,Sullivan~\cite{Sul71}. D.\,Sullivan
found obstructions to resolving singularities. These obstructions
appear to be elements of finite order. Hence D.\,Sullivan's
results imply that one can always resolve singularities ``with
multiplicity''. Our goal is to give an explicit construction for
such resolution of singularities.

It is convenient to work with pseudo-manifolds decomposed into
{\it simple cells}. A simple cell is a closed disk whose boundary
is decomposed into faces so that the decomposition is dual to a
piecewise linear triangulation of a sphere. Rigorous definitions
of a simple cell and a decomposition into simple cells can be
found in~\cite{Gai04}. Examples of decompositions into simple
cells are simplicial and cubic decompositions and cell
decompositions dual to piecewise linear triangulations of
manifolds. All manifolds and pseudo-manifolds are supposed to be
closed.

Suppose $Z$ is an oriented $n$-dimensional pseudo-manifold
decomposed into simple cells. The {\it problem of resolving
singularities} of a pseudo-manifold~$Z$ is to construct an
oriented piecewise linear manifold~$M$ and a mapping~$g:M\to Z$
such that off the $(n-2)$-dimensional skeleton of~$Z$ the
mapping~$g$ is a finite-fold covering.

{\it Main construction.} To each $(n-1)$-dimensional cell of the
decomposition~$Z$ we assign a label which is an element of a
certain finite set. The set of labels of $(n-1)$-dimensional
cells containing a cell~$F$ is called a {\it label} of~$F$ and is
denoted by~$c(F)$. A labeling is said to be {\it good} if

\begin{enumerate}
\item $|c(F)|=n-\dim F$ for any cell $F$ and

\item for any cell~$F$ and any $n$-dimensional
cell $G\supset F$ the labels of all cells $H$ such that $F\subset
H\subset G$ are pairwise distinct.
\end{enumerate}
\begin{lemma}
Suppose $Z$ is a pseudo-manifold decomposed into simple cells.
Then there is a pseudo-manifold $\overline{Z}$ decomposed into
simple cells such that $\overline{Z}$ admits a good labeling and
there is a mapping $p:\overline{Z}\to Z$ such that $p$ is a
covering off the codimension~$2$ skeleton and $p$ maps each cell
of $\overline{Z}$ isomorphically onto a cell of~$Z$.
\end{lemma}

Thus in the sequel we may assume that the decomposition~$Z$ is
endowed with a good labeling. Suppose $F$ is a cell of the
decomposition $Z$ such that $\dim F=k<n$. By $L_F$ we denote the
set of $n$-dimensional cells containing~$F$. Since $Z$ admits a
good labeling, cells $G\in L_F$ can be regularly colored in two
colours. (``Regularly'' means that any two cells possessing a
common facet are of distinct colours.) Besides, the numbers of
cells of both colours are equal to each other. By~$P_F$ we denote
the set of colour-reversing involutions on the set $L_F$. By~$P$
we denote the direct product of the sets~$P_F$ over all cells~$F$
of the decomposition~$Z$ such that $\dim F<n$.

By $U$ we denote the set of $(n+1)$-tuples $(F_0,F_1,\ldots,F_n)$,
where $F_0\supset F_1\supset\ldots\supset F_n$ are cells of~$Z$
and $\dim F_i=n-i$. We put, $V=U\times P\times\Bbb Z_2^n$.

We define involutions $\Phi^{\varepsilon}_j:V\to V$,
$\varepsilon=0,1$, $j=1,2,\ldots,n$, by
\begin{gather*}
\Phi_j^0(F_0,F_1,\ldots,F_n,(\Lambda_F),h)=
(F_0,F_1,\ldots,F_{j-1},F^*_j,F_{j+1}\ldots,F_n,(\Lambda_F),h),\\
\Phi_j^1(F_0,F_1,\ldots,F_n,(\Lambda_F),h)=
(\widetilde{F}_0,\widetilde{F}_1,\ldots,\widetilde{F}_n,(\widetilde{\Lambda}_F),h+e_j),
\end{gather*}
where $\Lambda_F\in P_F$, $h\in \Bbb Z_2^n$,
$(e_1,e_2,\ldots,e_n)$ is a basis of~$\Bbb Z_2^n$ and
\begin{enumerate}
\item $F^*_j$ is a unique cell such that $F_j^*\ne F_j$ and
$F_{j-1}\subset F^*_j\subset F_{j+1}$;
\item $\widetilde{F}_i=F_i$ if $i\geqslant j$;
$\widetilde{F}_0=\Lambda_{F_j}(F_0)$; if $0<i<j$, then
$\widetilde{F}_i$ is a unique cell such that
$\widetilde{F}_j\subset\widetilde{F}_i\subset\widetilde{F}_0$ and
$c(\widetilde{F}_i)=c(F_i)$;
\item if $c(F)\not\supset c(F_j)$, then
$\widetilde{\Lambda}_F=\Lambda_F$;
\item if $c(F)\supset c(F_j)$ and $G\in L_F$, then
$\widetilde{\Lambda}_F(G)=(\Lambda_{H_2}\circ\Lambda_F\circ\Lambda_{H_1})(G)$,
where $H_1$ is a unique cell such that~$F\subset H_1\subset G$ and
$c(H_1)=c(F_j)$ and $H_2$ is a unique cell such that~$F\subset
H_2\subset \Lambda_F(\Lambda_{H_1}(G))$ and~$c(H_2)=c(F_j)$.
\end{enumerate}
We put, $M=[0,1]^n\times V/\sim$, where $\sim$ is the equivalence
relation generated by the identifications
$$
(t_1,t_2,\ldots,t_n,\Phi_j^{\varepsilon}(v))\sim(t_1,t_2,\ldots,t_n,v)\text{
if }t_j=\varepsilon,\ \varepsilon=0,1,\ j=1,2,\ldots,n.
$$
It can be immediately checked that
$\Phi_j^0\Phi_k^1=\Phi_k^1\Phi_j^0$ if $j\ne k$ and
$\Phi_j^1\Phi_k^1=\Phi_k^1\Phi_j^1$. This implies that~$M$ is an
oriented piecewise linear manifold. A required mapping~$g:M\to Z$
is given by
$$
g(t_1,\ldots,t_n,F_0,\ldots,F_n,(\Lambda_F),h)=
\left(\prod_{i=1}^n(1-t_i)\right)b(F_0)+
\sum_{j=1}^n\left(t_j\prod_{i=j+1}^n(1-t_i)\right)b(F_j),
$$
where $b(F)$ is the barycenter of the cell~$F$.

Thus we obtain the following theorem.

\begin{theorem}
Main construction yields a mapping $g:M\to Z$ providing a
solution to the problem of resolving singularities of a
pseudo-manifold~$Z$. Suppose $X$ is a topological space and
$f:Z\to X$ is a singular cycle. Then the mapping $f\circ g:M\to
X$ realizes a multiple of the homology class of the cycle~$f$.
\end{theorem}

As an application of Theorem~1 we give a partial solution of the
following problem. Given a set of oriented simplicial spheres
$Y_1,Y_2,\ldots,Y_k$ construct an oriented simplicial manifold
whose set of links of vertices coincides up to an isomorphism
with the set $Y_1,Y_2,\ldots,Y_k$. Such manifold may exist only if
vertices of the triangulations~$Y_1,Y_2,\ldots,Y_k$ can be paired
off so that links of vertices in each pair are isomorphic to each
other with the isomorphism reversing the orientation. (If this
condition holds, the set $Y_1,Y_2,\ldots,Y_k$ is referred to be
{\it balanced}.)
\begin{theorem}
Suppose $Y_1,Y_2,\ldots,Y_k$ is a balanced set of oriented
simplicial spheres. There is an explicit construction yielding an
oriented simplicial manifold, whose set of links of vertices
coincides up to isomorphism with the set
$$
\underbrace{Y_1,\ldots,Y_1}_{r},\underbrace{Y_2,\ldots,Y_2}_{r},
\ldots,\underbrace{Y_k,\ldots,Y_k}_{r},K_1,K_2,\ldots,K_l,-K_1,-K_2,\ldots,-K_l,
$$
where $-K_i$ is a simplicial sphere $K_i$ with the opposite
orientation.
\end{theorem}

The author is grateful to V.\,M.\,Buchstaber for posing the
problems and permanent attention to his work.

\end{document}